\documentclass{amsart}
\usepackage{amsmath,amsfonts,amssymb}
\usepackage{epsfig}
\newcommand{\ie}{\emph{i.e.}}

\newcommand{\cf}{\emph{cf}}
\newcommand{\R}{\mathbb{R}}
\newcommand{\NM}{\mathbb{N}}
\newcommand{\s}{L}
\newcommand{\Smooth}{C}
\newcommand{\Dom}{\mathcal{D}}
\newcommand{\Hilbert}{\mathcal{H}}
\newcommand{\eps}{\varepsilon}
\newcommand{\axis}{r}
\newcommand{\domain}{\Omega}
\newcommand{\tube}{\mathcal{T}}
\newcommand{\nodal}{\mathcal{N}}
\newtheorem*{theo}{Theorem}
\newtheorem{conjecture}{}
\newtheorem{lem}{Lemma}
\newtheorem{prop}{Proposition}
\begin{document}
%
%
\title[Domains whose nodal line does not touch the boundary]
{Unbounded planar domains whose second nodal line does not touch the boundary}
\author{
Pedro Freitas and
David Krej\v{c}i\v{r}\'{\i}k
}
\thanks{
Partially supported by
FCT/POCTI/FEDER, Portugal,
and the Czech Academy of Sciences and its Grant Agency
within the projects IRP AV0Z10480505 and A100480501.
}
\subjclass{Primary 35P05; Secondary 35B05, 58J50}
\keywords{Dirichlet Laplacian, eigenfunctions,
nodal line, unbounded do\-mains}
\address{
Faculdade de Motricidade Humana {\tiny and}
Mathematical Physics Group of the University of Lisbon,
Complexo Interdisciplinar,
Av.~Prof.~Gama Pinto~2,
P-1649-003 Lisboa,
Portugal
}
\email{freitas@cii.fc.ul.pt}
\address{
Department of Theoretical Physics,
Nuclear Physics Institute, Academy of Sciences,
250\,68 \v{R}e\v{z} near Prague, Czech Republic
}
\email{krejcirik@ujf.cas.cz}
\date{14 November 2005;
\emph{accepted for publication in} Mathematical Research Letters}
\begin{abstract}
We show the existence of simply-connected
unbounded planar domains for which the
second nodal line of the Dirichlet Laplacian
does not touch the boundary.
\end{abstract}
\maketitle

%
\section{Introduction}
%
Consider the eigenvalue problem
\begin{equation}\label{maineq}
  \left\{
  \begin{aligned}
    -\Delta u &= \lambda u
    &\quad\mbox{in}\quad & \Omega \,,
    \\
    u &= 0
    &\quad\mbox{on}\quad & \partial \Omega \,,
  \end{aligned}
  \right.
\end{equation}
where~$\Omega$ is a domain (\ie\ open connected set) in~$\R^{2}$.
We interprete~(\ref{maineq}) in a weak sense as the eigenvalue problem
for the Dirichlet Laplacian $-\Delta_D^{\Omega}$
acting in the Hilbert space~$\s^2(\Omega)$,
and recall that $-\Delta_D^{\Omega}$ is
the non-negative self-adjoint operator
associated with the quadratic form
\begin{equation*}
  Q_D^\Omega[v] := \|\nabla v\|_{\s^2(\Omega)}^2 \,,
  \qquad
  v \in \Dom(Q_D^\Omega) := \Hilbert_0^1(\Omega) \,.
\end{equation*}

We denote by $\{\lambda_k(\Omega)\}_{k=1}^\infty$
the non-decreasing sequence of numbers corresponding
to the spectral problem of~$-\Delta_D^{\Omega}$
according to the Rayleigh-Ritz variational formula~\cite[Sec.~4.5]{Davies}.
Each~$\lambda_k(\Omega)$ represents either a discrete eigenvalue
or the threshold of the essential spectrum (if~$\Omega$ is not bounded).
All the eigenvalues below the essential spectrum
of the boundary-value problem~(\ref{maineq})
may be characterized by this variational principle.

The nodal line of a real eigenfunction $u$ of the problem~(\ref{maineq})
is defined by
\[
\nodal(u) = \overline{\left\{ x\in\Omega: u(x)=0\right\}} \,,
\]
and the connected components into which~$\Omega$ is divided by~$\nodal(u)$
are called the nodal domains of $u$.
By the Courant nodal domain theorem
an eigenfunction corresponding to the $k^\mathrm{th}$~eigenvalue
below the essential spectrum has at most~$k$
nodal domains (see~\cite{CH1} for the proof in the bounded case,
the generalization to the unbounded case being straightforward).
In particular, since the first eigenvalue below the essential spectrum
is always simple and the corresponding eigenfunction
can be chosen to be positive,
any eigenfunction corresponding to the second eigenvalue
below the essential spectrum will have exactly two nodal domains.

Apart from this result, very little is known regarding
the structure of the nodal lines,
but much work has been developed over the last three decades
around a conjecture of Payne's which states that a second
eigenfunction of the above problem cannot have a closed nodal line.
This conjecture is also quite often stated as follows:
\begin{conjecture}\label{conj}
The nodal line of any second eigenfunction of the
Laplacian intersects the boundary $\partial\Omega$ at exactly two points.
\end{conjecture}

The most general result obtained so far was given by Melas in 1992,
who showed that the above conjecture holds in the case of bounded
planar convex domains~\cite{Melas}. This followed a string of results obtained
under some symmetry restrictions by several authors (Payne himself
included) -- see, for instance,~\cite{Freitas1} and the references therein.

On the other hand, several counterexamples have also been presented,
of which the most significant is that in~\cite{H2ON} showing that the
result does not hold for multiply connected planar domains in general.
Other counterexamples have been given illustrating other ways in which
the conjecture may not hold. These include adding a
potential~(\cite{LinNi}) and the case of simply-connected
surfaces~(\cite{Freitas1}).

The purpose of this note is to give examples showing that
if one does not require the domain to be bounded, then the nodal line
need not touch the boundary even under the same assumptions that have
been previously used in the bounded case to prove the conjecture. More
precisely we will prove the following
\begin{theo} \label{t1}
    There exists a simply-connected unbounded planar
    domain~$\Omega$ which is convex and symmetric
    with respect to two orthogonal directions,
    and for which the nodal line of a second eigenfunction does not touch the
    boundary~$\partial\Omega$.

    This domain can be chosen as one of the following two types:
\begin{enumerate}
\item[(i)] the distance between the nodal line of a second
eigenfunction and the boundary $\partial\Omega$ is bounded away from zero,
but the spectrum is not purely discrete;

\item[(ii)] the spectrum consists only of discrete eigenvalues, but
the infimum of the distance between a point on the nodal line of a second
eigenfunction and the boundary $\partial\Omega$ is zero.
\end{enumerate}
\end{theo}

\begin{figure}[t]
\begin{tabular}{ll}
(i) & (ii) \\
\epsfig{file=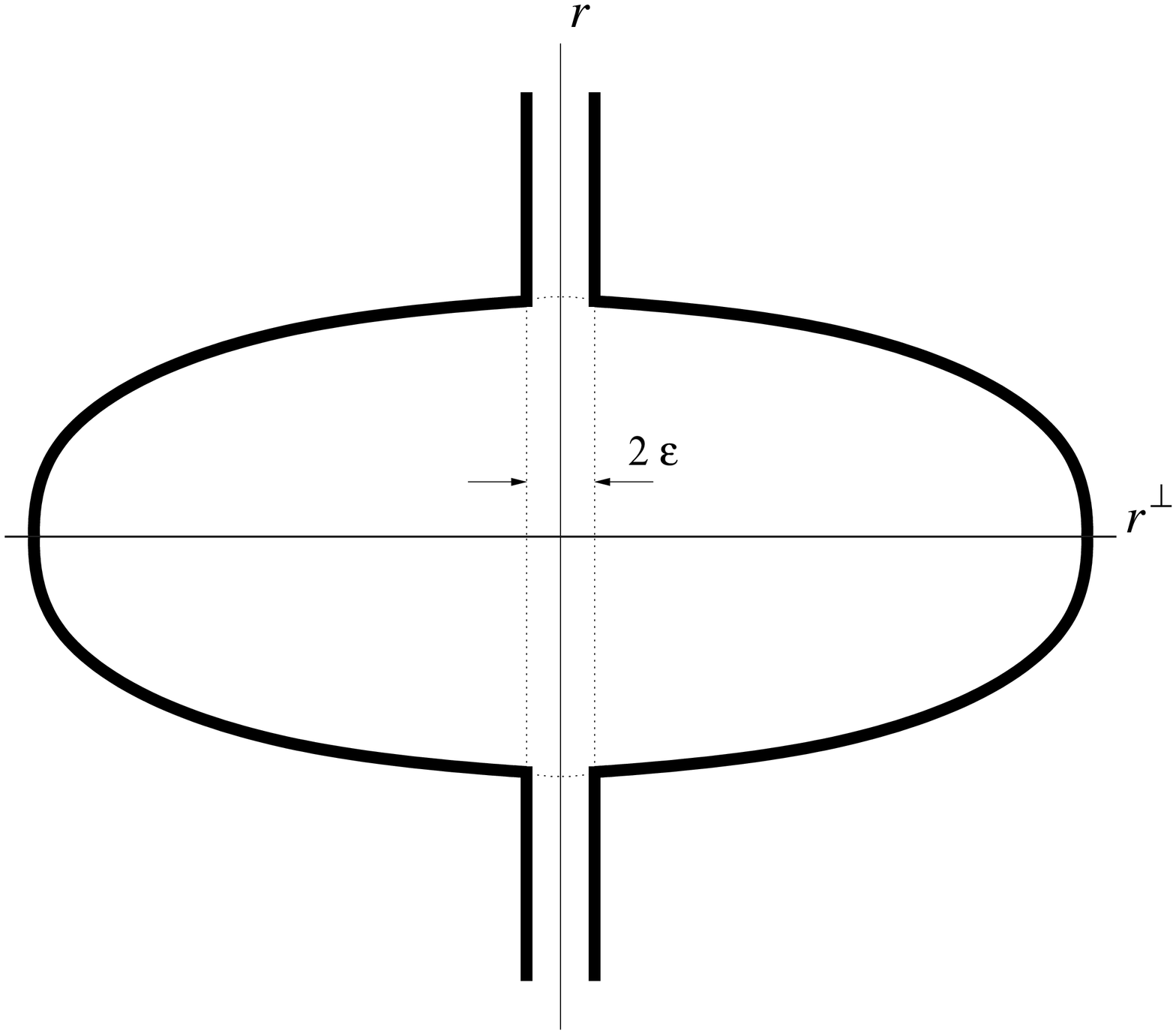,width=0.47\textwidth} &
\epsfig{file=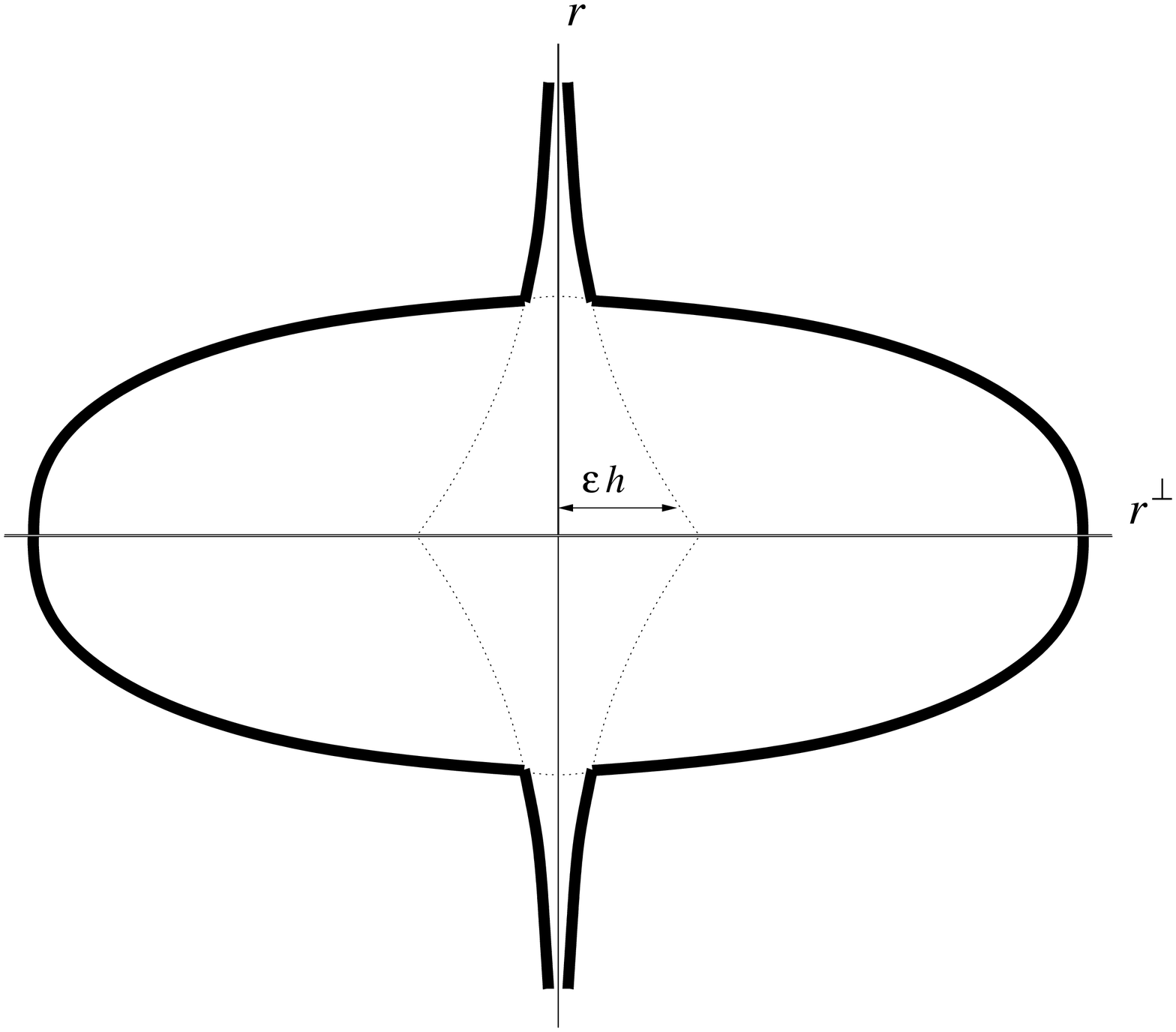,width=0.47\textwidth}
\end{tabular}
\caption{
Typical domains for which the nodal line of the second eigenfunction
does not touch the boundary.
}
\label{fig}
\end{figure}

The idea behind both examples is to start from
a bounded convex domain~$\domain_0$
which is invariant under reflections
through two orthogonal lines~$\axis$ and~$\axis^\bot$,
and which we will assume to be sufficiently long in the direction~$\axis^{\bot}$,
such that its second eigenvalue is simple
and any corresponding eigenfunction is antisymmetric with respect to~$\axis$.
In fact, its second nodal line will be
given by the closure of~$\domain_0 \cap \axis$.
We then append two sufficiently thin
semi-infinite strips to $\Omega_{0}$ in neighbourhoods of the points where
its second nodal line touches the boundary, in such a way that the
nodal line coincides with the axis~$\axis$
and thus stays within these strips without touching the boundary
-- see Figure~\ref{fig}.

In order to establish case~(i),
we will consider domains which are asymptotically cylindrical -- see the
classification of Euclidean domains
in~\cite[\S~49]{Glazman}, where these sets are called
quasi-cylindrical. This means that there will also exist essential
spectrum, and so it will be necessary to prove that the domain
does indeed possess a second discrete eigenvalue in this case.
In order for condition~(ii) to be satisfied,
we will need to consider what are referred to in~\cite{Glazman} as
quasi-bounded domains.
This means that the domains are asymptotically narrow and thus,
although the nodal line does not touch the boundary,
it does get asymptotically close to it.

It should be stressed here that while the nodal line
in both our examples does not touch the boundary,
it is not closed.

%

\section{The Proof}
%
Let~$\domain_0$ be a bounded open convex subset of~$\R^2$
which is simultaneously invariant under the reflection through
the coordinate axes
$
  \axis := \{0\} \times \R
$
and
$
  \axis^\bot := \R \times \{0\}
$,
\ie,
$$
  \forall (x_1,x_2) \in \R^2, \qquad
  (x_1,x_2)\in\domain_0
  \quad \Longrightarrow \quad
  \begin{cases}
    (x_1,-x_2)\in\domain_0 \,, \\
    (-x_1,x_2)\in\domain_0 \,.
  \end{cases}
$$
Of course, the first eigenvalue~$\lambda_1(\domain_0)>0$ is simple
and the corresponding eigenfunction can be chosen to be positive.
We assume that also the second eigenvalue~$\lambda_2(\domain_0)$
is simple and that the nodal line of the corresponding eigenfunction
is the closure of~$\domain_0\cap\axis$
(by~\cite{Jerison1,GJerison} and the symmetry,
these always happen if~$\domain_0$
is sufficiently long in the direction~$\axis^\bot$).

Let
$
  h:[0,+\infty)\to(0,1]
$
be a convex function.
For any~$\eps>0$,
we define an open tubular neighbourhood of the axis~$\axis$ by
$$
  \tube_\eps :=
  \left\{
  (x_1,x_2) \in \R^2 :\, |x_1| < \eps\,h(|x_2|)
  \right\} ,
$$
and introduce the unbounded open connected set
\begin{equation}
  \domain_\eps := \domain_0 \cup \tube_\eps
  \,.
\end{equation}
Note that~$\domain_\eps$ is invariant
both under the reflections through~$\axis$ and~$\axis^\bot$,
and convex both along~$\axis$ and~$\axis^\bot$.
It is also worth to notice that the boundary~$\partial\domain_\eps$
is necessarily at least of class~$\Smooth^{0,1}$,
\cf~\cite[Sec.~V.4.1]{Edmunds-Evans}.

Using the minimax principle
and a Dirichlet-Neumann bracketing argument,
it is easy to see that
\begin{equation}\label{sp.ess}
  \inf\sigma_\mathrm{ess}(-\Delta_D^{\Omega_\eps})
  \geq \pi^2/(2\eps)^2
\end{equation}
and that one can produce an arbitrary number
of eigenvalues below the essential spectrum by making~$\eps$ small enough.
Furthermore, if~$h$ tends to zero at infinity
then
$
  \sigma_\mathrm{ess}(-\Delta_D^{\Omega_\eps})
  = \varnothing
$
and the spectrum of~$-\Delta_D^{\Omega_\eps}$
consists of discrete eigenvalues only.
In any case, one has the following convergence result.
\begin{lem}\label{Lem.convergence}
$\forall k\in\NM\!\setminus\!\{0\}$, \quad
$
  {\displaystyle
  \lim_{\eps \to 0} \lambda_k(\domain_\eps) = \lambda_k(\domain_0)
  }
$.
\end{lem}
\begin{proof}
It follows from~\cite{Daners} that
$-\Delta_D^{\Omega_\eps}$ converges to $-\Delta_D^{\Omega_0}$
in the generalized sense of Kato~\cite{Kato},
which implies, in particular,
continuity of eigenvalues below the essential spectrum.
\end{proof}

In view of~(\ref{sp.ess}) and Lemma~\ref{Lem.convergence},
$\lambda_1(\domain_\eps)$ and~$\lambda_2(\domain_\eps)$
are discrete simple eigenvalues for all sufficiently small~$\eps$.
Since the eigenfunction corresponding
to~$\lambda_1(\domain_\eps)$ can be chosen to be positive,
the eigenfunction~$u_{2,\eps}$ corresponding to~$\lambda_2(\domain_\eps)$
has to change sign in~$\domain_\eps$.
We now prove a result which will give us immediately
the conclusions of the Theorem.
\begin{prop}
$\exists\eps_0>0$, $\forall\eps\in(0,\eps_0)$, \quad
$\nodal(u_{2,\eps}) = \axis$.
\end{prop}
\begin{proof}
Let~$\eps$ be so small that~$\lambda_2(\domain_\eps)$
is a simple discrete eigenvalue.
Due to the symmetry of~$\domain_\eps$,
the corresponding eigenfunction~$u_{2,\eps}$
must be symmetric or antisymmetric with respect to~$\axis$,
and symmetric or antisymmetric with respect to~$\axis^\bot$.
This observation and the Courant nodal domain theorem yield
that the nodal set~$\nodal(u_{2,\eps})$
is either the closure of~$\axis^\bot\cap\domain_0$,
the axis~$\axis$ or a closed loop.
The last possibility is excluded by mimicking
the argument given in~\cite{Damascelli} (see also~\cite{Payne2})
for bounded domains with the required symmetry and convexity.
We will exclude the first possibility by using the fact
that the second eigenvalue of a domain is the first eigenvalue
of any of the nodal subdomains.
Let us assume that there is a positive sequence $\{\eps_j\}_{j\in\NM}$,
converging to zero as $j\to\infty$,
such that
$
  \nodal(u_{2,\eps_j}) = \overline{\axis^\bot\cap\domain_0}
$
for all~$j\in\NM$.
Then
\begin{equation*}
  \lambda_2(\domain_{\eps_j})
  = \lambda_1\big(\domain_{\eps_j}\cap[\R\times(0,+\infty)]\big)
  \longrightarrow
  \lambda_1\big(\domain_{0}\cap[\R\times(0,+\infty)]\big)
  \qquad\mbox{as}\quad
  j\to\infty
\end{equation*}
by a convergence argument analogous to Lemma~\ref{Lem.convergence}.
On the other hand, we know that
\begin{equation*}
  \lambda_2(\domain_{\eps_j})
  \longrightarrow
  \lambda_2(\domain_{0})
  = \lambda_1\big(\domain_0\cap[(0,+\infty)\times\R]\big)
  \qquad\mbox{as}\quad
  j\to\infty
\end{equation*}
by Lemma~\ref{Lem.convergence}
and the assumption we have made about~$\domain_0$.
This implies that $\lambda_2(\domain_0)$ is degenerate
(there is one eigenfunction antisymmetric with respect to~$\axis$
and one eigenfunction antisymmetric with respect to~$\axis^\bot$),
a contradiction.
\end{proof}

It follows that the nodal line $\nodal(u_{2,\eps})$
does not touch the boundary of~$\domain_\eps$
for all sufficiently small~$\eps>0$.
Furthermore, if we choose~$h \equiv 1$
then the distance between the nodal line and the boundary is equal to~$\eps$,
which establishes part~(i) of the Theorem.
Part~(ii) follows by taking a function~$h$
which tends to zero at infinity.

%
\providecommand{\bysame}{\leavevmode\hbox to3em{\hrulefill}\thinspace}
\providecommand{\MR}{\relax\ifhmode\unskip\space\fi MR }
\providecommand{\MRhref}[2]{%
  \href{http://www.ams.org/mathscinet-getitem?mr=#1}{#2}
}
\providecommand{\href}[2]{#2}

\end{document}